\documentclass[reqno,12pt]{amsart}
\usepackage{amscd,amsfonts,mathrsfs,amsthm,amssymb, amsmath,mathtools}
\usepackage{upgreek}
\usepackage{csquotes}
\usepackage{enumerate}
\usepackage{bbm}
\usepackage{multicol}
\usepackage{stmaryrd,color}
\usepackage{array}
\usepackage{ae}
\usepackage{accents}
\usepackage{epsfig}
\usepackage[e]{esvect}
\usepackage[nobysame]{amsrefs}
\usepackage[margin=1in]{geometry}
\usepackage[toc,title,page]{appendix}
\usepackage{hyperref}
\hypersetup{colorlinks=true,linkcolor=blue,filecolor=blue,urlcolor=blue, citecolor=blue}

\let\oldtocsection=\tocsection
\let\oldtocsubsection=\tocsubsection
\let\oldtocsubsubsection=\tocsubsubsection
\renewcommand{\tocsection}[2]{\hspace{0em}\oldtocsection{#1}{#2}}
\renewcommand{\tocsubsection}[2]{\hspace{1em}\oldtocsubsection{#1}{#2}}
\renewcommand{\tocsubsubsection}[2]{\hspace{2em}\oldtocsubsubsection{#1}{#2}}

\newtheorem{theorem}{Theorem}

\newtheorem{lemma}{Lemma}
\newtheorem{remark}{Remark}

\makeatletter
\@namedef{subjclassname@2020}{\textup{2020} Mathematics Subject Classification}
\makeatother

\def\pproof#1{\@ifnextchar[\opargproof
{\opargproof[\it Proof of #1.]}}
\def\opargproof[#1]{\par\noindent {\bf #1 }}

\makeatother

\numberwithin{equation}{section}

\begin{document}

\date{\today}
\subjclass[2020]{Primary: 35B40, 35C20, 35K05, 35K65, 35R01, 47A11, 58J50. Secondary: 35P15, 47A10, 47B12, 58C40}
\thanks{The research project is implemented in the framework of H.F.R.I call ``Basic research Financing (Horizontal support of all Sciences)'' under the National Recovery and Resilience Plan ``Greece 2.0" funded by the European Union - NextGenerationEU (H.F.R.I. Project Number: 14758).}

\title[The spectrum of the Laplacian and the heat asymptotics]{The spectrum of the Laplacian on closed manifolds and the heat asymptotics near conical points}
\author[N. Roidos]{\textsc{Nikolaos Roidos}}
\address{Department of Mathematics, University of Patras, 26504 Rio Patras, Greece}
\email{roidos@math.upatras.gr}

\begin{abstract}
Let $\mathcal{M}$ be a smooth, closed and connected manifold of dimension $n\in\mathbb{N}$, endowed with a Riemannian metric $g$. Moreover, let $\mathcal{B}$ be an $(n+1)$-dimensional compact manifold with boundary equal to $\mathcal{M}$. Endow $\mathcal{B}$ with a Riemannian metric $h$ such that, in local coordinates $(x,y)\in [0,1)\times \mathcal{M}$ on the collar part of the boundary, it admits the warped product form $h=dx^{2}+x^{2}g(y)$. We consider the homogeneous heat equation on $(\mathcal{B},h)$ and find an arbitrary long asymptotic expansion of the solutions with respect to $x$ near $0$. It turns out that the spectrum of the Laplacian on $(\mathcal{M},g)$ determines explicitly the above asymptotic expansion and vice versa.
\end{abstract}

\maketitle

\section{Introduction}

Let $\mathcal{M}$ be a smooth, closed (i.e. compact without boundary) and connected manifold of dimension $n\in\mathbb{N}$, endowed with a Riemannian metric $g$. Denote $\mathbb{M}=(\mathcal{M},g)$ and let $\cdots<\lambda_{1}<\lambda_{0}=0$ be the spectrum of the Laplacian $\Delta_{\mathbb{M}}$ on $\mathbb{M}$. Moreover, let $\mathcal{B}$ be an $(n+1)$-dimensional compact manifold with boundary equal to $\mathcal{M}$, i.e. $\partial\mathcal{B}=\mathcal{M}$. Endow $\mathcal{B}$ with a Riemannian metric $h$ such that, in local coordinates $(x,y)\in [0,1)\times \mathcal{M}$ on the collar part of the boundary, it admits the warped product form $h=dx^{2}+x^{2}g(y)$. $\mathbb{B}=(\mathcal{B},h)$ is a manifold with an {\em isolated conical singularity} at $\{0\}\times \mathcal{M}$; $x$ is the distance from the singularity along the normal geodesics, see \cite{MW}.

Let $p\in(1,\infty)$, $s\in\mathbb{R}$,
 \begin{equation}\label{gammachoice}
\frac{n-3}{2}<\gamma<\min\Big\{\frac{n+1}{2},-1+\sqrt{\Big(\frac{n-1}{2}\Big)^{2}-\lambda_{1}}\Big\}
\end{equation}
 and let $\mathcal{H}_{p}^{s,\gamma}(\mathbb{B})$ be the Mellin-Sobolev space of order $s$ and weight $\gamma$ on $\mathcal{B}$, see Section 3. Let $\Delta_{\mathbb{B}}$ be the Laplacian on $\mathbb{B}$ and consider the homogeneous heat equation on $\mathbb{B}$, namely 
\begin{eqnarray}\label{heq1}
u'(t)&=&\Delta_{\mathbb{B}}u(t), \quad t>0,\\\label{heq2}
u(0)&=&u_{0}\in \mathcal{H}_{p}^{s,\gamma}(\mathbb{B}),
\end{eqnarray}
where $(\cdot)'=\partial_{t}$. Let $\omega$ be a fixed smooth non-negative real-valued function on $\mathbb{R}$, such that $\omega=1$ if $x\leq1/2$ and $\omega=0$ if $x\geq1$. Consider the space $\mathbb{C}_{\omega}=\{c\, \omega\, | \, c\in\mathbb{C}\}$ endowed with the norm $c\, \omega \mapsto |c|$. Then, the map 
\begin{equation}\label{deltadom}
\mathcal{H}_{p}^{s+2,\gamma+2}(\mathbb{B})\oplus\mathbb{C}_{\omega}\ni f \mapsto \Delta_{\mathbb{B}}f \in \mathcal{H}_{p}^{s,\gamma}(\mathbb{B})
\end{equation}
defines a closed extension of $\Delta_{\mathbb{B}}$ in $\mathcal{H}_{p}^{s,\gamma}(\mathbb{B})$, which we denote by $\underline{\Delta}_{\mathbb{B}}$, see Section 3. Moreover, $\underline{\Delta}_{\mathbb{B}}$ generates a holomorphic $C_{0}$-semigroup on $\mathcal{H}_{p}^{s,\gamma}(\mathbb{B})$, see the Appendix. As a consequence, see e.g. \cite[Corollary 3.7.21]{ABHN}, there exists a unique 
 \begin{equation}\label{solofheat}
 u\in C([0,+\infty);\mathcal{H}_{p}^{s,\gamma}(\mathbb{B}))\cap C^{\infty}((0,+\infty);\mathcal{H}_{p}^{s,\gamma}(\mathbb{B}))\cap C((0,+\infty);\mathcal{H}_{p}^{s+2,\gamma+2}(\mathbb{B})\oplus\mathbb{C}_{\omega})
 \end{equation}
solving \eqref{heq1}-\eqref{heq2}. For $m\in\mathbb{N}$ and $\gamma$ as in \eqref{gammachoice} denote
$$
\delta_{n,m}=\Big\{\begin{array}{lll} 1 &\text{if} &m=n \\ 0 & \text{if} & m\neq n \end{array} \quad \text{and} \quad J_{m}=\Big[\frac{n+1}{2}-\gamma-2m,\frac{n+1}{2}-\gamma-2(m-1)\Big),
$$
and let
\begin{equation}\label{specpowers}
\mu_{j}=\frac{n-1}{2}-\sqrt{\Big(\frac{n-1}{2}\Big)^{2}-\lambda_{j}}, \quad j\in\mathbb{N}.
\end{equation}
For the solution \eqref{solofheat} we show the following. 

\begin{theorem}\label{heateqexp}
Assume that $\lambda_{1}<(\frac{n-1}{2})^{2}-1$. Let $p\in(1,\infty)$, $s\in\mathbb{R}$,
\begin{equation}\label{gammaextra}
\max\Big\{\frac{n-3}{2},1-\sqrt{\Big(\frac{n-1}{2}\Big)^{2}-\lambda_{1}}\Big\}<\gamma<\min\Big\{\frac{n+1}{2},-1+\sqrt{\Big(\frac{n-1}{2}\Big)^{2}-\lambda_{1}}\Big\},
\end{equation}
and let $\underline{\Delta}_{\mathbb{B}}$ be the Laplacian \eqref{deltadom}. Then, the solution $u$ of the heat equation \eqref{heq1}-\eqref{heq2}, in addition to \eqref{solofheat}, satisfies
 \begin{equation}\label{uininftydom}
 u\in \bigcap_{k\in\mathbb{N}, k\geq2}C^{\infty}((0,+\infty);D(\underline{\Delta}_{\mathbb{B}}^{k})),
 \end{equation}
where for each $k\in\mathbb{N}$, $k\geq2$, for the domain of the $k$-th power we have
\begin{equation}\label{DDk}
D(\underline{\Delta}_{\mathbb{B}}^{k})=D(\underline{\Delta}_{\mathbb{B},\min}^{k})\oplus\mathbb{C}_{\omega}\oplus\Big(\bigoplus_{\nu=1}^{k-1}\big(\underline{\mathcal{F}}_{\nu\delta_{n,2}}^{(1-2\nu)\delta_{n,2}}\oplus\underline{\mathcal{F}}_{\nu+\delta_{n,1}}^{-2\nu}\big) +\sum_{m=2}^{k}\bigoplus_{\mu_{j}\in J_{m}} \bigoplus_{\nu=0}^{k-m}\underline{\mathcal{F}}_{m+\nu-2}^{-2\nu+\mu_{j}}\Big).
\end{equation} 
In the above equation the minimal domain satisfies
\begin{equation}\label{DDkminexp}
\mathcal{H}^{s+2k,\gamma+2k}_{p}(\mathbb{B}) \hookrightarrow D(\underline{\Delta}_{\mathbb{B},\min}^{k})\hookrightarrow \bigcap_{\varepsilon>0}\mathcal{H}^{s+2k,\gamma+2k-\varepsilon}_{p}(\mathbb{B})
\end{equation}
and, for $m\in\mathbb{N}_{0}=\mathbb{N}\cup\{0\}$ and $\rho\leq0$, by $\underline{\mathcal{F}}_{m}^{\rho}$ we denote a subspace of the finite-dimensional space
$$
\mathcal{F}_{m}^{\rho}=\Big\{\omega(x)x^{-\rho}\big(c_{0}(y)+c_{1}(y)\log(x)+\dots+c_{m}(y)\log^{m}(x)\big)\, |\, c_{j}\in C^{\infty}(\mathbb{M};\mathbb{C}), j\in\{0,\dots,m\}\Big\}.
$$
In particular, if $v\in D(\underline{\Delta}_{\mathbb{B},\min}^{k})$ and $s+2k>(n+1)/p$, then $v\in C(\mathbb{B})$ and for any $\varepsilon>0$, in local coordinates $(x,y)\in [0,1)\times \mathcal{M}$ on the collar part of the boundary, we have
\begin{equation}\label{MelSobemb}
|v(x,y)|\leq Cx^{\gamma+2k-\frac{n+1}{2}-\varepsilon}\|v\|_{\mathcal{H}^{s+2k,\gamma+2k-\varepsilon}_{p}(\mathbb{B})}, 
\end{equation}
for some $C>0$ depending on $\mathbb{M}$, $n$ and $p$.
\end{theorem}

Note that the above theorem is free of the assumption when $n\geq3$. The regularity \eqref{uininftydom} for solutions of the heat equation on conic manifolds was first considered in \cite[Theorem 4.3]{Roidos}, without providing a description of the domain of the $k$-th power of the Laplacian. The description \eqref{DDk}, which is the main result of this article, has the following consequence. Let $p\in(1,\infty)$, $s\in\mathbb{R}$, $k\in\mathbb{N}$, $k\geq2$, satisfying $s+2k>(n+1)/p$, $\gamma$ be as in \eqref{gammaextra} and let $u_{0}\in \mathcal{H}_{p}^{s,\gamma}(\mathbb{B})$. Then, according to Theorem \ref{heateqexp}, there exist sequences $\{a_{\nu}\}_{\nu\in\mathbb{N}}$, $\{b_{\nu+\delta_{n,1}}\}_{\nu\in\mathbb{N}}$, $\{c_{m,j,m+\nu-3}\}_{\nu\in\mathbb{N}}$, $m,j\in\mathbb{N}$, $m\geq2$, depending on $\mathbb{M}$, $n$, $\gamma$ and $u_{0}$, with 
$$
a_{\eta}, b_{\eta}, c_{m,\eta}\in\Big\{\omega(x)\sum_{j=0}^{\eta}\theta_{j}(t,y)\log^{j}(x) \, |\, \theta_{j}(t,\cdot)\in C^{\infty}(\mathbb{M};\mathbb{C}), t>0, j\in\{0,\dots,\eta\}\Big\},
$$ 
such that, for the solution $u$ of the heat equation \eqref{heq1}-\eqref{heq2} we have
\begin{eqnarray}\nonumber
u(t,x,y)&=&\sum_{\nu=1}^{k-1}\big(\delta_{n,2}a_{\nu}(t,x,y)x^{2\nu-1}+b_{\nu+\delta_{n,1}}(t,x,y)x^{2\nu}\big)\\\nonumber
&&\hspace{-54pt}+\sum_{m=2}^{k}\,\sum_{(\frac{n-1}{2})^{2}-(\gamma+2m-1)^{2}\leq \lambda_{j}<(\frac{n-1}{2})^{2}-(\gamma+2m-3)^{2}} x^{-\frac{n-1}{2}+\sqrt{\big(\frac{n-1}{2}\big)^{2}-\lambda_{j}}} \sum_{\nu=m-2}^{k-2} c_{m,j,\nu}(t,x,y)x^{2(\nu-m+2)} \\\label{uexp}
&&\hspace{-54pt}+\, \omega(x)(c(t)+v(t,x,y)),
\end{eqnarray}
$t>0$, $(x,y)\in[0,1)\times\mathcal{M}$, where $c\in C^{\infty}((0,\infty);\mathbb{C})$, $v\in C^{\infty}((0,\infty);C(\mathbb{B}))$ and for any $\varepsilon>0$ we have 
$$
|v(t,x,y)|\leq Lx^{\gamma+2k-\frac{n+1}{2}-\varepsilon},
$$
for some constant $L$ depending only on $\mathbb{M}$, $n$, $p$, $s$, $\gamma$, $k$, $\varepsilon$, $t$ and $v$. 

Based on \eqref{uexp}, we see that the spectrum of $\Delta_{\mathbb{M}}$ determines explicitly the asymptotic behaviour of $u$ near the conical tip. On the other hand, if we know the asymptotic behaviour of $u$ near $0$, by \eqref{uexp} we can recover the full spectrum of $\Delta_{\mathbb{M}}$. 

As an example, let us choose $\mathbb{M}$ to be the $n$-dimensional sphere $\mathbb{S}^{n}$, $n\geq2$. In this case the spectrum of the Laplacian is $\lambda_{j}=-j(j+n-1)$, $j\in\mathbb{N}_{0}$. Therefore, $\mu_{j}=-j$, $j\in\mathbb{N}$, so that \eqref{uexp} provides the usual Taylor expansion for $u$ near $x=0$.

\section{Holomorphic semigroups and their associated Cauchy problems}

Let $X$ be a complex Banach space and $A:D(A)\rightarrow X$ be a closed and densely defined linear operator in $X$. Assume that there exist some $\theta\in(\pi/2,\pi)$ and $K\geq1$ such that
$$
S_{\theta}=\{0\}\cup\{z\in\mathbb{C}\backslash\{0\}\, | \, |\arg{z}|\leq \theta\} \subset \rho(-A) 
$$
and
$$
\|(A+z)^{-1}\|_{\mathcal{L}(X)}\leq \frac{K}{1+|z|} \quad \text{when} \quad z\in S_{\theta},
$$
where by $\mathcal{L}(\cdot)$ we denote the space of bounded linear operators. Then, $-A$ generates a bounded holomorphic $C_{0}$-semigroup $\{e^{-tA}\}_{t\geq0}$ on $X$, given by the inverse Laplace transform representation
\begin{equation}\label{etA}
e^{-tA}=\frac{1}{2\pi i}\int_{\Gamma_{\theta}}e^{tz}(A+z)^{-1}dz, \quad t>0,
\end{equation}
where $\Gamma_{\theta}$ is the upwards oriented boundary of $S_{\theta}$, see e.g. \cite[Theorem 13.30]{Neerven} or \cite[Theorem 3.3.1]{Tanabe}. As a consequence, for any $w_{0}\in X$ there exists a unique 
$$
w\in C^{\infty}((0,+\infty);X)\cap C([0,+\infty);X)\cap C((0,+\infty);D(A))
$$
satisfying 
$$
w'(t)+Aw(t)=0, \quad t>0, \quad w(0)=w_{0};
$$
it is given by $w(t)=e^{-tA}w_{0}$, $t\geq0$, see e.g. \cite[Corollary 3.7.21]{ABHN}. The above solution also satisfies 
\begin{equation}\label{extrareg}
w\in \bigcap_{k\in\mathbb{N}}C^{\infty}((0,+\infty);D(A^{k})),
\end{equation}
see e.g. \cite[Remark 3.7.20]{ABHN} or \cite[Proposition 2.1.1]{Lunardi}, where for each $k\in\mathbb{N}$ the $k$-th power $A^{k}$ of $A$ is defined as usual by
$$
D(A^{m+1})=\{u\in D(A^{m})\, |\, Au \in D(A^{m})\}, \quad m\in\mathbb{N}.
$$
The regularity \eqref{extrareg} can be seen as follows. Since for any $t>0$ the integral 
$$
\int_{\Gamma_{\theta}}e^{tz}A(A+z)^{-1}w_{0}dz
$$
converges absolutely, by \eqref{etA} we have that $e^{-tA}w_{0}\in D(A)$ and
\begin{eqnarray*}
Ae^{-tA}w_{0}&=&\frac{1}{2\pi i}\int_{\Gamma_{\theta}}e^{tz}A(A+z)^{-1}w_{0}dz\\
&=&\Big(\frac{1}{2\pi i}\int_{\Gamma_{\theta}}e^{tz}dz\Big)w_{0}-\frac{1}{2\pi i}\int_{\Gamma_{\theta}}e^{tz}z(A+z)^{-1}w_{0}dz\\
&=&-\frac{1}{2\pi i}\int_{\Gamma_{\theta}}e^{tz}z(A+z)^{-1}w_{0}dz,
\end{eqnarray*}
where we have used Cauchy's theorem. Since the integral 
$$
\int_{\Gamma_{\theta}}e^{tz}zA(A+z)^{-1}w_{0}dz
$$
converges absolutely, we have that $e^{-tA}w_{0}\in D(A^{2})$ and
\begin{eqnarray*}
A^{2}e^{-tA}w_{0}&=&-\frac{1}{2\pi i}\int_{\Gamma_{\theta}}e^{tz}zA(A+z)^{-1}w_{0}dz\\
&=&\Big(-\frac{1}{2\pi i}\int_{\Gamma_{\theta}}e^{tz}zdz\Big)w_{0}+\frac{1}{2\pi i}\int_{\Gamma_{\theta}}e^{tz}z^{2}(A+z)^{-1}w_{0}dz\\
&=&\frac{1}{2\pi i}\int_{\Gamma_{\theta}}e^{tz}z^{2}(A+z)^{-1}w_{0}dz.
\end{eqnarray*}
After $m$ steps we deduce that $e^{-tA}w_{0}\in D(A^{m})$ and 
$$
A^{m}e^{-tA}w_{0}=\frac{(-1)^{m}}{2\pi i}\int_{\Gamma_{\theta}}e^{tz}z^{m}(A+z)^{-1}w_{0}dz.
$$
Hence 
\begin{equation}\label{wexpAm}
w(t)=e^{-tA}w_{0}=\frac{(-1)^{m}}{2\pi i}\int_{\Gamma_{\theta}}e^{tz}z^{m}(A+z)^{-1}A^{-m}w_{0}dz, \quad t>0.
\end{equation}
Therefore, if $k\in \mathbb{N}$, since the graph norm of $D(A^{k})$ is equivalent to $\|A^{k}\cdot\|_{X}$, by taking $m\geq k$, for $\tau>0$ sufficiently small we have
\begin{eqnarray*}
\lefteqn{\|w(t+\tau)-w(t)\|_{D(A^{k})}}\\
&=&\|\frac{(-1)^{m}}{2\pi i}\int_{\Gamma_{\theta}}(e^{\tau z}-1)e^{tz}z^{m}(A+z)^{-1}A^{-m}w_{0}dz\|_{D(A^{k})}\\
&\leq&C_{0}\|\frac{(-1)^{m}}{2\pi i}\int_{\Gamma_{\theta}}(e^{\tau z}-1)e^{tz}z^{m}(A+z)^{-1}A^{k-m}w_{0}dz\|_{X},
\end{eqnarray*}
for certain $C_{0}>0$. By the dominated convergence theorem, the last term in the above inequality tends to zero as $\tau\rightarrow 0$. Thus,
\begin{equation}\label{wCAk}
w\in C((0,+\infty);D(A^{k})).
\end{equation}
Furthermore, 
\begin{eqnarray*}
\lefteqn{\|\frac{1}{\tau}(w(t+\tau)-w(t))+Aw(t)\|_{D(A^{k})}}\\
&=&\|\frac{(-1)^{m}}{2\pi i}\int_{\Gamma_{\theta}}\frac{e^{\tau z}-1}{\tau}e^{tz}z^{m}(A+z)^{-1}A^{-m}w_{0}dz\\
&&+\frac{(-1)^{m+1}}{2\pi i}\int_{\Gamma_{\theta}}e^{tz}z^{m+1}(A+z)^{-1}A^{-m}w_{0}dz\|_{D(A^{k})}\\
&=&\|\frac{(-1)^{m}}{2\pi i}\int_{\Gamma_{\theta}}\Big(\frac{e^{\tau z}-1}{\tau}-z\Big)e^{tz}z^{m}(A+z)^{-1}A^{-m}w_{0}dz\|_{D(A^{k})}\\
&\leq&C_{0}\|\frac{(-1)^{m}}{2\pi i}\int_{\Gamma_{\theta}}\Big(\frac{e^{\tau z}-1}{\tau}-z\Big)e^{tz}z^{m}(A+z)^{-1}A^{k-m}w_{0}dz\|_{X}.
\end{eqnarray*}
By the dominated convergence theorem, the last term in the above inequality tends to zero as $\tau\rightarrow 0$. Hence, due to \eqref{wCAk}, we get
$$
w\in C^{1}((0,+\infty);D(A^{k})) \quad \text{and} \quad w'(t)=-Aw(t), \quad t>0.
$$
By \eqref{wexpAm} we have
$$
w'(t)=\frac{(-1)^{m+1}}{2\pi i}\int_{\Gamma_{\theta}}e^{tz}z^{m}(A+z)^{-1}A^{1-m}w_{0}dz.
$$
Therefore, by taking $m\geq k+1$, we estimate
\begin{eqnarray*}
\lefteqn{\|\frac{1}{\tau}(w'(t+\tau)-w'(t))-A^{2}w(t)\|_{D(A^{k})}}\\
&=&\|\frac{(-1)^{m+1}}{2\pi i}\int_{\Gamma_{\theta}}\frac{e^{\tau z}-1}{\tau}e^{tz}z^{m}(A+z)^{-1}A^{1-m}w_{0}dz\\
&&-\frac{(-1)^{m+1}}{2\pi i}\int_{\Gamma_{\theta}}e^{tz}z^{m+1}(A+z)^{-1}A^{1-m}w_{0}dz\|_{D(A^{k})}\\
&=&\|\frac{(-1)^{m+1}}{2\pi i}\int_{\Gamma_{\theta}}\Big(\frac{e^{\tau z}-1}{\tau}-z\Big)e^{tz}z^{m}(A+z)^{-1}A^{1-m}w_{0}dz\|_{D(A^{k})}\\
&\leq&C_{0}\|\frac{(-1)^{m+1}}{2\pi i}\int_{\Gamma_{\theta}}\Big(\frac{e^{\tau z}-1}{\tau}-z\Big)e^{tz}z^{m}(A+z)^{-1}A^{k+1-m}w_{0}dz\|_{X}.
\end{eqnarray*}
By the dominated convergence theorem, the last term in the above inequality tends to zero as $\tau\rightarrow 0$. Thus, due to \eqref{wCAk}, we get
$$
w\in C^{2}((0,+\infty);D(A^{k})) \quad \text{and} \quad w''(t)=A^{2}w(t), \quad t>0.
$$
After $\nu$ steps, we deduce that 
$$
w\in C^{\nu}((0,+\infty);D(A^{k})) \quad \text{and} \quad \Big(\frac{d}{dt}\Big)^{\nu}w(t)=(-1)^{\nu}A^{\nu}w(t), \quad t>0,
$$
which confirms \eqref{extrareg}.

\section{Cone differential operators}

In this section we recall some basic facts about cone differential operators and their associated pseudo-differential theory, which is called {\em cone calculus}. For details we refer to \cite{CSS1}, \cite{Gill}, \cite{GKM2}, \cite{GKM1}, \cite{GM}, \cite{Lesch}, \cite{SS1}, \cite{SS2}, \cite{Schulze1}, \cite{Schulze2} and \cite{Schulze3}. A $\mu$-th order differential operator $A$ with smooth coefficients in the interior $\mathbb{B}^{\circ}$ of $\mathbb{B}$ is called {\em cone differential operator} (with $x$-independent coefficients near the boundary) if in local coordinates $(x,y)\in [0,1)\times \mathcal{M}$ on the collar part it is expressed as
\begin{equation}\label{coneop}
A=x^{-\mu}\sum_{j=0}^{\mu}\alpha_{j}(-x\partial_{x})^{j}, \quad \text{where} \quad \alpha_{j} \in \text{Diff}^{\mu-j}(\mathbb{M}), \quad j\in\{0,\dots,\mu\}.
\end{equation}
Denote by $C_{c}^{\infty}(\cdot)$ the space of smooth compactly supported functions, $H_{p}^{s}(\cdot)$, $p\in(1,\infty)$, $s\in\mathbb{R}$, the usual {\em Bessel potential space} and let $H_{p,loc}^{s}(\cdot)$ be the space of {\em locally $H_{p}^{s}(\cdot)$-functions}. Cone differential operators act naturally on scales of Mellin-Sobolev spaces. Let $\kappa_{j}:\Omega_{j}\subseteq\mathcal{M} \rightarrow\mathbb{R}^{n}$, $j\in\{1,\dots,\ell\}$, $\ell\in\mathbb{N}$, be a covering of $\mathcal{M}$ by coordinate charts and let $\{\phi_{j}\}_{j\in\{1,\dots,\ell\}}$ be a subordinate partition of unity. For any $p\in(1,\infty)$ and $s,\gamma\in\mathbb{R}$ let $\mathcal{H}^{s,\gamma}_p(\mathbb{B})$ be the space of all distributions $u$ on $\mathbb{B}^{\circ}$ such that 
$$
\|u\|_{\mathcal{H}^{s,\gamma}_p(\mathbb{B})}=\sum_{j=1}^{\ell}\|G_{\gamma}(1\otimes \kappa_{j})_{\ast}(\omega\phi_{j} u)\|_{H^{s}_{p}(\mathbb{R}^{n+1})}+\|(1-\omega)u\|_{H^{s}_{p}(\mathbb{B})}
$$
is defined and finite, where 
$$
G_{\gamma}: C_{c}^{\infty}(\mathbb{R}_{+}\times\mathbb{R}^{n})\rightarrow C_{c}^{\infty}(\mathbb{R}^{n+1}) \quad \mbox{is defined by} \quad u(x,y)\mapsto e^{(\gamma-\frac{n+1}{2})x}u(e^{-x},y)
$$
and $\ast$ refers to the push-forward of distributions. The space $\mathcal{H}^{s,\gamma}_{p}(\mathbb{B})$, called {\em (weighted) Mellin-Sobolev space}, is independent of the choice of the cut-off function $\omega$, the covering $\{\kappa_{j}\}_{j\in\{1,\dots,\ell\}}$ and the partition $\{\phi_{j}\}_{j\in\{1,\dots,\ell\}}$. If in particular $s\in \mathbb{N}_{0}$, then equivalently, $\mathcal{H}^{s,\gamma}_{p}(\mathbb{B})$, $p\in(1,\infty)$, $\gamma\in\mathbb{R}$, is the space of all functions $u$ in $H^s_{p,loc}(\mathbb{B}^\circ)$ such that near the boundary we have
$$
x^{\frac{n+1}{2}-\gamma}(x\partial_{x})^{k}\partial_{y}^{\alpha}(\omega u) \in L^{p}\big([0,1)\times\mathcal{M},\sqrt{\det[g]}\frac{dx}{x}dy\big),\quad k+|\alpha|\leq s.
$$
Moreover, if $A$ is as in \eqref{coneop}, then for any $\eta\in\mathbb{R}$ the map
$$
A:\mathcal{H}^{\eta+\mu,\gamma+\mu}_p(\mathbb{B})\rightarrow \mathcal{H}^{\eta,\gamma}_p(\mathbb{B})
$$
is well defined and bounded. 

Beyond the usual {\em homogeneous principal symbol} $\sigma_{\psi}(A)\in C^{\infty}((T^{\ast}\mathbb{B}^{\circ})\backslash\{0\} )$ of $A$ in \eqref{coneop}, its {\em principal rescaled symbol} is defined by 
$$
\widetilde{\sigma}_{\psi}(A)=\sum_{j=0}^{\mu}\sigma_{\psi}(\alpha_{j})(y,\xi)(-i\tau)^{j} \in C^{\infty}((T^{\ast}\mathbb{M}\times \mathbb{R})\backslash\{0\}),
$$
where $ (x,y,\tau,\xi)$ are local coordinates in $T^{\ast}([0,1)\times\mathbb{M})$. The operator $A$ is called {\em $\mathbb{B}$-elliptic} if both $\sigma_{\psi}(A)$ and $\widetilde{\sigma}_{\psi}(A)$ are pointwise invertible.

Assume now that $A$ in \eqref{coneop} is $\mathbb{B}$-elliptic and regard it as an unbounded operator in $\mathcal{H}^{s,\gamma}_{p}(\mathbb{B})$, $p\in(1,\infty)$, $s,\gamma\in\mathbb{R}$, with domain $C_{c}^{\infty}(\mathbb{B}^{\circ})$. The domain of the minimal extension (i.e. the closure) $\underline{A}_{\min}$ of $A$ is given by
$$
D(\underline{A}_{\min})=\Big\{u\in\bigcap_{\varepsilon>0} \mathcal{H}^{s+\mu,\gamma+\mu-\varepsilon}_{p}(\mathbb{B}) \, |\, Au\in \mathcal{H}^{s,\gamma}_{p}(\mathbb{B})\Big\}.
$$
By \cite[Proposition 2.3]{SS2} we have 
\begin{equation}\label{DAminemb}
 \mathcal{H}^{s+\mu,\gamma+\mu}_{p}(\mathbb{B})\hookrightarrow D(\underline{A}_{\min}) \hookrightarrow \mathcal{H}^{s+\mu,\gamma+\mu-\varepsilon}_{p}(\mathbb{B})
\end{equation}
for all $\varepsilon>0$ and $D(\underline{A}_{\min})=\mathcal{H}^{s+\mu,\gamma+\mu}_{p}(\mathbb{B})$ if and only if $\sigma_{M}(A)(\cdot)$ is invertible on the line $\{z\in\mathbb{C} \, | \, \mathrm{Re}(z)=\frac{n+1}{2}-\gamma-\mu\}$, where the family of differential operators
$$
\mathbb{C}\ni z \mapsto \sigma_{M}(A)(z)=\sum_{j=0}^{\mu}\alpha_{j}z^{j} \in \mathcal{L}(H_{2}^{\mu}(\mathbb{M}),H_{2}^{0}(\mathbb{M}))
$$
is called the {\em conormal symbol} of $A$. Note that, if $B$ is a cone differential operator of order $\nu$, then
\begin{equation}\label{symbAB}
\sigma_{M}(AB)(z)=\sigma_{M}(A)(z+\nu)\sigma_{M}(B)(z), \quad z\in\mathbb{C},
\end{equation}
see e.g. \cite[Proposition 2.4]{Gill} or \cite[Lemma 1.1.13 (2)]{Lesch} or \cite[(2.13)]{SS2}.

On the other hand, the domain of the maximal extension $\underline{A}_{\max}$ of $A$, defined as usual by 
$$
D(\underline{A}_{\max})=\{u\in \mathcal{H}^{s,\gamma}_{p}(\mathbb{B})\, |\, Au\in \mathcal{H}^{s,\gamma}_{p}(\mathbb{B})\},
$$
is given by 
\begin{equation}\label{maxdomexp}
D(\underline{A}_{\max})=D(\underline{A}_{\min})\oplus\bigoplus_{\frac{n+1}{2}-\gamma-\mu\leq \rho < \frac{n+1}{2}-\gamma}\mathcal{E}_{\rho}.
\end{equation}
Here $\rho\in [\frac{n+1}{2}-\gamma-\mu,\frac{n+1}{2}-\gamma)$ runs over all poles of $( \sigma_{M}(A)(\cdot))^{-1}$ and, for each such $\rho$, $\mathcal{E}_{\rho}$ is a $(p,s)$-independent finite-dimensional space consisting of linear combinations of functions that vanish on $\mathcal{B}\backslash([0,1)\times\mathcal{M})$ and, in local coordinates $(x,y)\in[0,1)\times\mathcal{M}$ on the collar part they are of the form $\omega(x)c(y)x^{-\rho}\log^{\nu}(x)$, $\nu\in (0,\dots,\eta_{\rho}-1)$, where $c \in C^{\infty}(\mathbb{M};\mathbb{C})$ and $\eta_{\rho}\in \mathbb{N}$ is the order of the pole $\rho$, see e.g. \cite[(3.4)]{CSS1} or \cite[Theorem 4.7]{GKM2} or \cite[(6.1)]{GKM1}. As a consequence, there are several closed extensions of $A$ in $\mathcal{H}^{s,\gamma}_{p}(\mathbb{B})$, also called realizations; each one is obtained from \eqref{maxdomexp} after choosing a subspace of 
$$
\bigoplus_{\frac{n+1}{2}-\gamma-\mu\leq \rho < \frac{n+1}{2}-\gamma}\mathcal{E}_{\rho}.
$$

We focus now on the Laplacian $\Delta_{\mathbb{B}}$. In local coordinates $(x,y)\in[0,1)\times\mathcal{M}$ on the collar part it takes the form
\begin{equation}\label{laplexp}
\Delta_{\mathbb{B}}=\partial_{x}^{2}+\frac{n}{x}\partial_{x}+\frac{1}{x^{2}}\Delta_{\mathbb{M}}=\frac{1}{x^{2}}\big((x\partial_{x})^{2}+(n-1)(x\partial_{x})+\Delta_{\mathbb{M}}\big).
\end{equation}
Hence, $\Delta_{\mathbb{B}}$ is a second order $\mathbb{B}$-elliptic cone differential operator whose conormal symbol is given by 
\begin{equation}
\sigma_{\Delta}(z)=z^{2}-(n-1)z+\Delta_{\mathbb{M}}, \quad z\in\mathbb{C}.
\end{equation}
The poles of $(\sigma_{\Delta}(\cdot))^{-1}$ are identified with the set
\begin{equation}\label{polesD}
q_{j}^{\pm}=\frac{n-1}{2}\pm\sqrt{\Big(\frac{n-1}{2}\Big)^{2}-\lambda_{j}}, \quad j\in\mathbb{N}_{0};
\end{equation}
we have 
$$
\dots<q_{1}^{-}<q_{0}^{-}=0\leq q_{0}^{+}=n-1<q_{1}^{+}<\cdots.
$$
If we choose $(n-3)/2<\gamma<(n+1)/2$, then $0\in [(n+1)/2-\gamma-2,(n+1)/2-\gamma)$ so that $\mathbb{C}_{\omega}$ belongs to the maximal domain of $\Delta_{\mathbb{B}}$ in $\mathcal{H}^{s,\gamma}_{p}(\mathbb{B})$, $p\in(1,\infty)$, $s\in\mathbb{R}$. Therefore, the map \eqref{deltadom} indeed defines a closed extension of $\Delta_{\mathbb{B}}$. The fact that $\underline{\Delta}_{\mathbb{B}}$ generates a holomorphic semigroup on $\mathcal{H}^{s,\gamma}_{p}(\mathbb{B})$ is due to the further restriction $\gamma<-1+((n-1)/2)^{2}-\lambda_{1})^{1/2}$, which makes no $\mathcal{E}_{q_{j}^{-}}$, $j\in\mathbb{N}$, contribution to the maximal domain of $\Delta_{\mathbb{B}}$, see the Appendix for details.

\section{Proof of Theorem \ref{heateqexp}}

Setting $u(t)=e^{\xi t}w(t)$ in \eqref{heq1}-\eqref{heq2}, $\xi>0$, we obtain
\begin{eqnarray*}
w'(t)+(\xi-\Delta_{\mathbb{B}})w(t)&=&0, \quad t>0,\\
w(0)&=&u_{0}.
\end{eqnarray*}
Hence, by \eqref{extrareg} we get that $u$, in addition to \eqref{solofheat}, satisfies
$$
u\in \bigcap_{k\in\mathbb{N}}C^{\infty}((0,+\infty);D((\xi-\underline{\Delta}_{\mathbb{B}})^{k})).
$$
Then \eqref{uininftydom} follows from the fact that the norm $\|(\xi-\underline{\Delta}_{\mathbb{B}})^{k}\cdot \|_{\mathcal{H}^{s,\gamma}_{p}(\mathbb{B})}$ is equivalent to the graph norms of $(\xi-\underline{\Delta}_{\mathbb{B}})^{k}$ and $\underline{\Delta}_{\mathbb{B}}^{k}$.

For deriving \eqref{DDk} we start with the following continuous embedding
\begin{equation}\label{DDmemb}
D(\underline{\Delta}_{\mathbb{B}}^{m+1})\hookrightarrow D(\underline{\Delta}_{\mathbb{B}}^{m}), \quad m\in\mathbb{N}.
\end{equation}
Note that $\Delta^{k}$ is a $\mathbb{B}$-elliptic cone differential operator of order $2k$ whose conormal symbol, due to \eqref{symbAB}, is given by 
\begin{equation}\label{consymnDk}
\sigma_{\Delta^{k}}(z)=\sigma_{\Delta}(z)\cdots\sigma_{\Delta}(z+2(k-1)), \quad z\in\mathbb{C}.
\end{equation}

Let $Q_{k}$ be the poles of $(\sigma_{\Delta^{k}}(\cdot))^{-1}$, which are determined by \eqref{polesD}, \eqref{consymnDk} and \cite[(2.10)]{RS3}; note that $q_{j}^{-}=\mu_{j}$, $j\in\mathbb{N}$. Recall that, for any $a\in\mathbb{R}$, $m\in\mathbb{N}_{0}$ and $c_{0},\dots,c_{m}\in C(\mathbb{M};\mathbb{C})$, we have
$$
\omega(x)x^{a}\sum_{\nu=0}^{m}c_{\nu}(y)\log^{\nu}(x)\in \mathcal{H}_{p}^{s,\gamma}(\mathbb{B}) \quad \text{if and only if} \quad a>\gamma-\frac{n+1}{2}.
$$
We start by examining the poles \eqref{polesD} of $(\sigma_{\Delta}(\cdot))^{-1}$ in $J_{1}$. Since the weight $\gamma$ satisfies \eqref{gammachoice}, we always have $0\in J_{1}$ and 
\begin{equation}\label{qjminus}
q_{j}^{-}<\frac{n+1}{2}-\gamma-2, \quad j\in\mathbb{N}.
\end{equation}
Moreover, the extra restriction \eqref{gammaextra} on the weight, i.e. the fact that
$$
1-\sqrt{\Big(\frac{n-1}{2}\Big)^{2}-\lambda_{1}}<\gamma,
$$
implies $q_{j}^{+}>\frac{n+1}{2}-\gamma$, $j\in\mathbb{N}$; this is a key point in our analysis since only $q_{j}^{-}$ shifted by even integers contribute now to $D(\underline{\Delta}_{\mathbb{B}}^{k})$. We conclude that $q_{j}^{\pm}\notin J_{1}$ when $j\in\mathbb{N}$. Finally, for the pole $q_{0}^{+}=n-1$ we have the following cases:\\
(i) If $n=1$, then $q_{0}^{+}=0\in J_{1}$; in this case $0$ is a double pole of $(\sigma_{\Delta}(\cdot))^{-1}$.\\
(ii) If $n=2$, then we can either have $q_{0}^{+}\in J_{1}$ or $q_{0}^{+}\notin J_{1}$, depending on the choice of $\gamma$.\\
(iii) If $n\geq3$, then we always have $q_{0}^{+}>\frac{n+1}{2}-\gamma$; and hence $q_{0}^{+}\notin J_{1}$.

We proceed now by induction. In order to make the induction step more clear, we study the cases of $k=2$, $3$ and $4$ separately.

{\em Case $k=2$}. By \eqref{maxdomexp} we have
$$
D(\underline{\Delta}_{\mathbb{B}}^{2})=D(\underline{\Delta}_{\mathbb{B},\min}^{2})\oplus\bigoplus_{q\in Q_{2}\cap(J_{1}\cup J_{2})}\underline{\mathcal{F}}_{1+\delta_{n,1}}^{q}.
$$
Due to \eqref{DDmemb} with $m=1$ and \eqref{DDkminexp} with $k=2$, we have
$$
D(\underline{\Delta}_{\mathbb{B}}^{2})=D(\underline{\Delta}_{\mathbb{B},\min}^{2})\oplus\mathbb{C}_{\omega}\oplus\bigoplus_{q\in Q_{2}\cap J_{2}}\underline{\mathcal{F}}_{1+\delta_{n,1}}^{q}.
$$
By \eqref{consymnDk} with $k=2$ we have that $-2\in Q_{2}\cap J_{2}$, and in the case of $n=2$, $q_{0}^{+}-2=-1\in Q_{2}\cap J_{2}$. Hence, due to \eqref{qjminus} we have
\begin{equation}\label{DD2}
D(\underline{\Delta}_{\mathbb{B}}^{2})=D(\underline{\Delta}_{\mathbb{B},\min}^{2})\oplus\mathbb{C}_{\omega}\oplus\Big(\underline{\mathcal{F}}_{\delta_{n,2}}^{-\delta_{n,2}}\oplus\underline{\mathcal{F}}_{1+\delta_{n,1}}^{-2}+\bigoplus_{q_{j}^{-}\in J_{2}}\underline{\mathcal{F}}_{0}^{q_{j}^{-}}\Big),
\end{equation}
where ``$+$" appears since we may have $q_{j}^{-}\in\{-2,-1\}$ for some $q_{j}^{-}\in J_{2}$. Moreover, the lower index $1+\delta_{n,1}$ in $\underline{\mathcal{F}}_{1+\delta_{n,1}}^{-2}$ is due to the fact that, when $n=1$ the pole $-2$ can be of order $3$. Finally, note that the term 
$$
\bigoplus_{q_{j}^{-}\in J_{2}}\underline{\mathcal{F}}_{0}^{q_{j}^{-}}
$$ 
in \eqref{DD2} is omitted when $q_{1}^{-}< \frac{n+1}{2}-\gamma-4$. 

{\em Case $k=3$}. By \eqref{maxdomexp} we have
$$
D(\underline{\Delta}_{\mathbb{B}}^{3})=D(\underline{\Delta}_{\mathbb{B},\min}^{3})\oplus\bigoplus_{q\in Q_{3}\cap(J_{1}\cup J_{2}\cup J_{3})}\underline{\mathcal{F}}_{2+\delta_{n,1}}^{q}.
$$
Due to \eqref{DDmemb} with $m=2$, \eqref{DDkminexp} with $k=3$ and \eqref{DD2}, we have
\begin{equation}\label{DDA3}
D(\underline{\Delta}_{\mathbb{B}}^{3})=D(\underline{\Delta}_{\mathbb{B},\min}^{3})\oplus\mathbb{C}_{\omega}\oplus\Big(\underline{\mathcal{F}}_{\delta_{n,2}}^{-\delta_{n,2}}\oplus\underline{\mathcal{F}}_{1+\delta_{n,1}}^{-2}+\bigoplus_{q_{j}^{-}\in J_{2}}\underline{\mathcal{F}}_{0}^{q_{j}^{-}}\Big)\oplus\bigoplus_{q\in Q_{3}\cap J_{3}}\underline{\mathcal{F}}_{2+\delta_{n,1}}^{q}.
\end{equation}
By \eqref{consymnDk} with $k=3$ we have that $Q_{3}\cap J_{3}$ is the union of the following sets:\\
(i) $\{q_{j}^{-}\, |\, q_{j}^{-}\in J_{3}\}\cup\{-2+q_{j}^{-}\, |\, q_{j}^{-}\in J_{2}\}$.\\
(ii) $\{-4,-3\}$, which appears due to $0,1\in J_{1}$ when they are shifted by $-4$. \\
Hence,
$$
\bigoplus_{q\in Q_{3}\cap J_{3}}\underline{\mathcal{F}}_{2+\delta_{n,1}}^{q}=\underline{\mathcal{F}}_{2\delta_{n,2}}^{-3\delta_{n,2}}\oplus\underline{\mathcal{F}}_{2+\delta_{n,1}}^{-4}+\bigoplus_{q_{j}^{-}\in J_{2}}\underline{\mathcal{F}}_{1}^{-2+q_{j}^{-}}+\bigoplus_{q_{j}^{-}\in J_{3}}\underline{\mathcal{F}}_{1}^{q_{j}^{-}}.
$$
Therefore, by \eqref{DDA3} we get
\begin{eqnarray}\nonumber
D(\underline{\Delta}_{\mathbb{B}}^{3})&=&D(\underline{\Delta}_{\mathbb{B},\min}^{3})\oplus\mathbb{C}_{\omega}\oplus\Big(\underline{\mathcal{F}}_{\delta_{n,2}}^{-\delta_{n,2}}\oplus\underline{\mathcal{F}}_{1+\delta_{n,1}}^{-2}\oplus \underline{\mathcal{F}}_{2\delta_{n,2}}^{-3\delta_{n,2}}\oplus\underline{\mathcal{F}}_{2+\delta_{n,1}}^{-4}\\\label{DD3}
&&+\bigoplus_{q_{j}^{-}\in J_{2}}\big(\underline{\mathcal{F}}_{0}^{q_{j}^{-}}\oplus\underline{\mathcal{F}}_{1}^{-2+q_{j}^{-}}\big)+\bigoplus_{q_{j}^{-}\in J_{3}}\underline{\mathcal{F}}_{1}^{q_{j}^{-}}\Big).
\end{eqnarray} 

{\em Case $k=4$}. By \eqref{maxdomexp} we have
$$
D(\underline{\Delta}_{\mathbb{B}}^{4})=D(\underline{\Delta}_{\mathbb{B},\min}^{4})\oplus\bigoplus_{q\in Q_{4}\cap(J_{1}\cup J_{2}\cup J_{3}\cup J_{4})}\underline{\mathcal{F}}_{3+\delta_{n,1}}^{q}.
$$
Due to \eqref{DDmemb} with $m=3$, \eqref{DDkminexp} with $k=4$ and \eqref{DD3}, we have
\begin{eqnarray}\nonumber
D(\underline{\Delta}_{\mathbb{B}}^{4})&=&D(\underline{\Delta}_{\mathbb{B},\min}^{4})\oplus\mathbb{C}_{\omega}\oplus\Big(\underline{\mathcal{F}}_{\delta_{n,2}}^{-\delta_{n,2}}\oplus\underline{\mathcal{F}}_{1+\delta_{n,1}}^{-2}\oplus \underline{\mathcal{F}}_{2\delta_{n,2}}^{-3\delta_{n,2}}\oplus\underline{\mathcal{F}}_{2+\delta_{n,1}}^{-4}\\\label{DADA4}
&&+\bigoplus_{q_{j}^{-}\in J_{2}}\big(\underline{\mathcal{F}}_{0}^{q_{j}^{-}}\oplus\underline{\mathcal{F}}_{1}^{-2+q_{j}^{-}}\big)+\bigoplus_{q_{j}^{-}\in J_{3}}\underline{\mathcal{F}}_{1}^{q_{j}^{-}}\Big)\oplus\bigoplus_{q\in Q_{4}\cap J_{4}}\underline{\mathcal{F}}_{3+\delta_{n,1}}^{q}.
\end{eqnarray} 
By \eqref{consymnDk} with $k=4$ we have that $Q_{4}\cap J_{4}$ is the union of the following sets:\\
(i) $\{q_{j}^{-}\, |\, q_{j}^{-}\in J_{4}\}\cup\{-2+q_{j}^{-}\, |\, q_{j}^{-}\in J_{3}\}\cup\{-4+q_{j}^{-}\, |\, q_{j}^{-}\in J_{2}\}$.\\
(ii) $\{-6,-5\}$, which appears due to $0,1\in J_{1}$ when they are shifted by $-6$. \\
Hence,
$$
\bigoplus_{q\in Q_{4}\cap J_{4}}\underline{\mathcal{F}}_{3+\delta_{n,1}}^{q}=\underline{\mathcal{F}}_{3\delta_{n,2}}^{-5\delta_{n,2}}\oplus\underline{\mathcal{F}}_{3+\delta_{n,1}}^{-6}+\bigoplus_{q_{j}^{-}\in J_{2}}\underline{\mathcal{F}}_{2}^{-4+q_{j}^{-}}+\bigoplus_{q_{j}^{-}\in J_{3}}\underline{\mathcal{F}}_{2}^{-2+q_{j}^{-}}+\bigoplus_{q_{j}^{-}\in J_{4}}\underline{\mathcal{F}}_{2}^{q_{j}^{-}}.
$$
Therefore, by \eqref{DADA4} we get
\begin{eqnarray*}
D(\underline{\Delta}_{\mathbb{B}}^{4})&=&D(\underline{\Delta}_{\mathbb{B},\min}^{4})\oplus\mathbb{C}_{\omega}\oplus\Big(\underline{\mathcal{F}}_{\delta_{n,2}}^{-\delta_{n,2}}\oplus\underline{\mathcal{F}}_{1+\delta_{n,1}}^{-2}\oplus \underline{\mathcal{F}}_{2\delta_{n,2}}^{-3\delta_{n,2}}\oplus\underline{\mathcal{F}}_{2+\delta_{n,1}}^{-4}\oplus\underline{\mathcal{F}}_{3\delta_{n,2}}^{-5\delta_{n,2}}\oplus\underline{\mathcal{F}}_{3+\delta_{n,1}}^{-6}\\
&&+\bigoplus_{q_{j}^{-}\in J_{2}}\big(\underline{\mathcal{F}}_{0}^{q_{j}^{-}}\oplus\underline{\mathcal{F}}_{1}^{-2+q_{j}^{-}}\oplus\underline{\mathcal{F}}_{2}^{-4+q_{j}^{-}}\big)+\bigoplus_{q_{j}^{-}\in J_{3}}\big(\underline{\mathcal{F}}_{1}^{q_{j}^{-}}\oplus\underline{\mathcal{F}}_{2}^{-2+q_{j}^{-}}\big)+\bigoplus_{q_{j}^{-}\in J_{4}}\underline{\mathcal{F}}_{2}^{q_{j}^{-}}\Big).
\end{eqnarray*} 

{\em From $k$ to $k+1$}. Assume that \eqref{DDk} holds for some $k\in\mathbb{N}$, $k\geq2$. By \eqref{maxdomexp} we have
$$
D(\underline{\Delta}_{\mathbb{B}}^{k+1})=D(\underline{\Delta}_{\mathbb{B},\min}^{k+1})\oplus\bigoplus_{q\in Q_{k+1}\cap(J_{1}\cup \cdots \cup J_{k+1})}\underline{\mathcal{F}}_{k+\delta_{n,1}}^{q}.
$$
Due to \eqref{DDmemb} with $m=k$, \eqref{DDkminexp} with $k$ replaced by $k+1$ and \eqref{DDk}, we have
\begin{eqnarray}\nonumber
D(\underline{\Delta}_{\mathbb{B}}^{k+1})&=&D(\underline{\Delta}_{\mathbb{B},\min}^{k+1})\oplus\mathbb{C}_{\omega}\oplus\Big(\bigoplus_{\nu=1}^{k-1}\big(\underline{\mathcal{F}}_{\nu\delta_{n,2}}^{(1-2\nu)\delta_{n,2}}\oplus\underline{\mathcal{F}}_{\nu+\delta_{n,1}}^{-2\nu}\big)\\\label{DADA1k1}
&&+\sum_{m=2}^{k}\bigoplus_{q_{j}^{-}\in J_{m}} \bigoplus_{\nu=0}^{k-m}\underline{\mathcal{F}}_{m+\nu-2}^{-2\nu+q_{j}^{-}}\Big)\oplus\bigoplus_{q\in Q_{k+1}\cap J_{k+1}}\underline{\mathcal{F}}_{k+\delta_{n,1}}^{q}.
\end{eqnarray}
By \eqref{consymnDk} with $k$ replaced by $k+1$ we have that $Q_{k+1}\cap J_{k+1}$ is the union of the following sets:\\
(i) 
$$
\bigcup_{m=2}^{k+1}\{-2(k+1-m)+q_{j}^{-} \, |\, q_{j}^{-}\in J_{m}\}.
$$
(ii) $\{-2k,1-2k\}$, which appears due to $0,1\in J_{1}$ when they are shifted by $-2k$. \\
Hence,
$$
\bigoplus_{q\in Q_{k+1}\cap J_{k+1}}\underline{\mathcal{F}}_{k+\delta_{n,1}}^{q}=\underline{\mathcal{F}}_{k\delta_{n,2}}^{(1-2k)\delta_{n,2}}\oplus\underline{\mathcal{F}}_{k+\delta_{n,1}}^{-2k}+\sum_{m=2}^{k+1}\bigoplus_{q_{j}^{-}\in J_{m}}\underline{\mathcal{F}}_{k-1}^{-2(k+1-m)+q_{j}^{-}}.
$$
Therefore, by \eqref{DADA1k1} we obtain \eqref{DDk} with $k$ replaced by $k+1$. The embedding \eqref{DDkminexp} follows by \eqref{DAminemb}. Finally, \eqref{MelSobemb} follows by \eqref{DDkminexp} and the Mellin-Sobolev embedding \cite[Corollary 2.9]{RS2}.

\section{Appendix}

Taking into account \cite[Corollary 3.7.17]{ABHN}, the fact $\underline{\Delta}_{\mathbb{B}}$ generates a holomorphic semigroup on $\mathcal{H}^{s,\gamma}_{p}(\mathbb{B})$, for any $p\in(1,\infty)$, $s\in\mathbb{R}$ and $\gamma$ as in \eqref{gammachoice}, follows by \cite[Theorem 4.2]{RS1} or \cite[Theorem 6.7]{SS1}. However, the proofs of these two theorems are quite complicated, since they concern {\em boundedness of the imaginary powers} and {\em boundedness of the $H^{\infty}$-functional calculus} for the Laplacian, i.e. properties beyond holomorphic semigroup generation. For this reason, and for the completeness of this article, we include here an elementary and independent proof of the holomorphic semigroup generation when $s\geq0$. The proof uses only basic facts of the cone calculus such as the descriptions of the maximal domain and the domain of the adjoint for a $\mathbb{B}$-elliptic cone differential operator.
 
Denote $\mathbb{Y}=([0,\infty)\times\mathcal{M}, dx^{2}+x^{2}g(y))$. Note that the Laplacian $\Delta_{\mathbb{Y}}$ on $\mathbb{Y}$ is given by \eqref{laplexp}. For any $p\in(1,\infty)$ and $s,\gamma\in\mathbb{R}$ let 
$$
\mathcal{K}^{s,\gamma}_{p}(\mathbb{Y})=\big\{u\,\, |\,\, \omega u\in \mathcal{H}^{s,\gamma}_{p}(\mathbb{B}) \quad \text{and} \quad (1-\omega)u\in H^{s}_{p,cone}((\mathbb{R}\times\mathcal{M},dx^{2}+g(y))\big\} ;
$$
for the definition of the $H^{s}_{p,cone}$-space we refer to \cite[Section 2.1.1]{SS1}. Similarly to $\underline{\Delta}_{\mathbb{B}}$, when $\gamma$ satisfies \eqref{gammachoice} the map
\begin{equation}\label{deltahat}
\mathcal{K}^{2,\gamma+2}_{p}(\mathbb{Y})\oplus\mathbb{C}_{\omega}\ni u\mapsto \Delta_{\mathbb{Y}}u\in \mathcal{K}^{0,\gamma}_{p}(\mathbb{Y})
\end{equation}
defines a closed extension of $\Delta_{\mathbb{Y}}: C_{c}^{\infty}(\mathbb{Y}^{\circ})\rightarrow \mathcal{K}^{0,\gamma}_{p}(\mathbb{Y})$, which we denote by $\underline{\Delta}_{\mathbb{Y}}$.

\begin{lemma}\label{nospecminus}
Let $p\in(1,\infty)$ and $\gamma$ satisfying \eqref{gammachoice}. Then, the spectrum of $\underline{\Delta}_{\mathbb{Y}}$ is contained in $(-\infty,0]$.
\end{lemma}
\begin{proof}
Let $\lambda\in \mathbb{C}\backslash(-\infty,0]$ and let $\{e_{j}\}_{j\in\mathbb{N}_{0}}$ be a complete ortonormal system in $L^{2}(\mathbb{M})$ such that $e_{j}\in C^{\infty}(\mathbb{M})$ and $\Delta_{\mathbb{M}}e_{j}=\lambda_{k_{j}}e_{j}$, $j\in\mathbb{N}_{0}$, for some monotonically non-decreasing onto map $k_{(\cdot)}:\mathbb{N}_{0}\rightarrow \mathbb{N}_{0}$. By the $\mathbb{B}$-ellipticity of $\Delta_{\mathbb{B}}$, the invertibility of $\lambda-\underline{\Delta}_{\mathbb{Y}}$ is independent of $p$, see e.g. \cite[Proposition 3.3]{SS1}. Hence, it suffices to consider the case of $p=2$. Assume that there exists a $u\in \mathcal{K}^{2,\gamma+2}_{2}(\mathbb{Y})\oplus\mathbb{C}_{\omega}$ such that $\underline{\Delta}_{\mathbb{Y}}u=\lambda u$. Then
$$
u=\sum_{j=0}^{\infty}a_{j}(x)e_{j}(y) \quad \text{and} \quad x^{2}a''_{j}(x)+nxa'_{j}(x)+(\lambda_{k_{j}}-\lambda x^{2})a_{j}(x)=0, \quad j\in\mathbb{N}_{0}.
$$
Hence, if we denote 
$$
r_{j}=\sqrt{\Big(\frac{n-1}{2}\Big)^{2}-\lambda_{k_{j}}}, \quad j\in\mathbb{N}_{0},
$$ 
then
\begin{equation}\label{ajexp}
a_{j}(x)=x^{\frac{1-n}{2}}\big(c_{1,j}J_{r_{j}}(-i x\sqrt{\lambda})+c_{2,j}Y_{r_{j}}(-i x\sqrt{\lambda})\big), \quad j\in\mathbb{N}_{0},
\end{equation}
for certain coefficients $c_{1,j}, c_{2,j}\in\mathbb{C}$, $j\in\mathbb{N}_{0}$, where $J_{\nu}$ and $Y_{\nu}$, $\nu\in\mathbb{R}$, stand for the Bessel functions of order $\nu$, of first and second kind respectively. Note that we always have 
\begin{equation}\label{arg0pi}
-\pi<\arg(-ix\sqrt{\lambda})<0.
\end{equation}
From the asymptotics of Bessel functions for large arguments (see e.g. \cite[Section 7.21, eq. (1)]{Wat}), we have that
\begin{eqnarray*}
\lefteqn{a_{j}(x)=x^{-\frac{n}{2}}\sqrt{\frac{2i}{\pi\sqrt{\lambda}}}\bigg[c_{1,j}\cos\Big(-ix\sqrt{\lambda}-\frac{\pi}{2}r_{j}-\frac{\pi}{4}\Big)}\\
&&+c_{2,j}\sin\Big(-ix\sqrt{\lambda}-\frac{\pi}{2}r_{j}-\frac{\pi}{4}\Big)+e^{|\mathrm{Re}(x\sqrt{\lambda})|}\mathcal{O}(x^{-1})\bigg], \quad j\in\mathbb{N}_{0}.
\end{eqnarray*}
Due to $u\in \mathcal{K}^{0,\gamma}_{2}(\mathbb{Y})$, taking into account \eqref{arg0pi} we find that $c_{1,j}=ic_{2,j}$, $j\in\mathbb{N}_{0}$. On the other hand, from the asymptotics of Bessel functions for small arguments (see e.g. \cite[Section 3.1, eq. (8) and Section 3.53, eq. (1)]{Wat} and \cite[(8.444)]{GR}), we have
\begin{eqnarray}\nonumber
\lefteqn{\hspace{-55pt}a_{j}(x)=x^{\frac{1-n}{2}}\bigg[\frac{c_{1,j}(\frac{-ix\sqrt{\lambda}}{2})^{r_{j}}}{\Gamma(1+r_{j})}-c_{2,j}\bigg(\frac{\Gamma(r_{j})}{\pi}\Big(\frac{2}{-ix\sqrt{\lambda}}\Big)^{r_{j}}}\\\label{asyatzero}
&&-\frac{\cot(\pi r_{j})}{\Gamma(1+r_{j})}\Big(\frac{-ix\sqrt{\lambda}}{2}\Big)^{r_{j}}\bigg)\bigg]+\text{smaller terms}, \quad j\in\mathbb{N}_{0},
\end{eqnarray}
if $(n,j)\neq(1,0)$ and
$$
a_{0}(x)=c_{1,0}+c_{2,0}\frac{2}{\pi}\Big(\ln\Big(\frac{-i x\sqrt{\lambda}}{2}\Big)+\ell\Big)+\text{smaller terms} \quad \text{if} \quad n=1,
$$
where $\ell$ stands for the Euler-Mascheroni constant. Due to the choice $\gamma>(n-3)/2$, we have that 
$$
\omega(x)\ln(x), \omega(x)x^{\frac{1-n}{2}-\sqrt{\big(\frac{n-1}{2}\big)^{2}-\lambda_{k_{j}}}}\notin \mathcal{K}^{2,\gamma+2}_{2}(\mathbb{Y}), \quad j\in \mathbb{N}_{0},
$$
so that $c_{2,j}=0$, $j\in\mathbb{N}_{0}$. We conclude that $\mathrm{Ker}(\lambda-\underline{\Delta}_{\mathbb{Y}})=\{0\}$.

The inner product $\langle\cdot,\cdot\rangle_{0,0}$ in $\mathcal{K}^{0,0}_{2}(\mathbb{Y})$ induces an identification of the dual space of $\mathcal{K}^{0,\gamma}_{2}(\mathbb{Y})$ with $\mathcal{K}^{0,-\gamma}_{2}(\mathbb{Y})$. The domain of the adjoint $\underline{\Delta}_{\mathbb{Y}}^{\ast}$ of $\underline{\Delta}_{\mathbb{Y}}$, is defined as usual by
\begin{eqnarray*}
\lefteqn{D(\underline{\Delta}_{\mathbb{Y}}^{\ast})=\big\{v\in \mathcal{K}_{2}^{0,-\gamma}(\mathbb{Y}) \,\, | \,\, \text{there exists a} \, \, w\in \mathcal{K}_{2}^{0,-\gamma}(\mathbb{Y})}\\
&&\hspace{40pt}\text{such that} \,\, \langle v, \Delta_{\mathbb{Y}}\eta\rangle_{0,0}=\langle w,\eta\rangle_{0,0} \,\, \text{for all} \,\, \eta\in D(\underline{\Delta}_{\mathbb{Y}})\big\}.
\end{eqnarray*}
In addition, $D(\underline{\Delta}_{\mathbb{Y}}^{\ast})\hookrightarrow D(\underline{\Delta}_{\mathbb{Y},\max}^{\ast})$ with the maximal domain expressed as
$$
D(\underline{\Delta}_{\mathbb{Y},\max}^{\ast})=\mathcal{K}_{2}^{2,2-\gamma}(\mathbb{Y})\oplus\bigoplus_{\frac{n-3}{2}+\gamma\leq q_{j}^{\pm} < \frac{n+1}{2}+\gamma}\mathcal{E}_{q_{j}^{\pm}},
$$
where $q_{j}^{\pm}$ are given by \eqref{polesD} and
$$
\mathcal{E}_{q_{j}^{\pm}}=\{(x,y)\mapsto \omega(x)x^{-q_{j}^{\pm}}\vartheta(y) \, |\, \text{$\vartheta$ is a $\lambda_{j}$-eigenfunction of $\Delta_{\mathbb{M}}$}\}, \quad j\in\mathbb{N},
$$
$$
\mathcal{E}_{q_{0}^{\pm}}=\bigg\{\begin{array}{lll} \{ (x,y)\mapsto \omega(x)c\, x^{q_{0}^{\pm}} \, |\, c \in\mathbb{C}\} & \text{if} & n>1 \\ 
\{(x,y)\mapsto \omega(x)(c+\underline{c} \log(x)) \, | \, c, \underline{c}\in \mathbb{C}\} & \text{if} & n=1, \end{array}
$$
are the associated asymptotics spaces (see e.g. \cite[Section 5.2]{SS2}). Assume that there exists a $v\in D(\underline{\Delta}_{\mathbb{Y}}^{\ast})$ such that $\underline{\Delta}_{\mathbb{Y}}^{\ast}v=\lambda v$. Then $v=\sum_{j=0}^{\infty}b_{j}(x)e_{j}(y)$ and $b_{j}$ satisfy the same equation as $a_{j}$. Hence, $b_{j}$ are given by \eqref{ajexp}, for some constants $\underline{c}_{1,j}$, $\underline{c}_{2,j}$ instead of $c_{1,j}$, $c_{2,j}$. Similarly to $c_{1,j}$, $c_{2,j}$, we get $\underline{c}_{1,j}=i\underline{c}_{2,j}$, $j\in\mathbb{N}_{0}$. Furthermore, due to the choice $\gamma<(n+1)/2$, we have that 
$$
\omega(x)x^{\frac{1-n}{2}-\sqrt{\big(\frac{n-1}{2}\big)^{2}-\lambda_{k_{j}}}}\notin \mathcal{K}^{2,-\gamma+2}_{2}(\mathbb{Y}), \quad j\in \mathbb{N}_{0},
$$
and
$$
\omega(x)\ln(x)\notin \mathcal{K}^{2,-\gamma+2}_{2}(\mathbb{Y}) \quad \text{when} \quad n=1.
$$
Moreover, due to the choice $\gamma<-1+\sqrt{\big(\frac{n-1}{2}\big)^{2}-\lambda_{1}}$, we have $q_{k_{j}}^{+}>\frac{n+1}{2}+\gamma$, $j\in \mathbb{N}$, so that
$$
\omega(x)x^{\frac{1-n}{2}-\sqrt{\big(\frac{n-1}{2}\big)^{2}-\lambda_{k_{j}}}}=\omega(x)x^{-q_{k_{j}}^{+}}\notin\bigoplus_{\frac{n-3}{2}+\gamma\leq q_{m}^{\pm} < \frac{n+1}{2}+\gamma}\mathcal{E}_{q_{m}^{\pm}}, \quad j\in \mathbb{N}.
$$
Now if $\omega(x)\log(x)$, when $n=1$, or $\omega(x)x^{1-n}$, when $n>1$, belongs to $D(\underline{\Delta}_{\mathbb{Y}}^{\ast})$, then by setting $\eta=\omega(x)$ and $v=\omega(x)\log(x)$, when $n=1$, or $v=\omega(x)x^{1-n}$, when $n>1$, in $\langle v, \Delta_{\mathbb{Y}}\eta\rangle_{0,0}= \langle \Delta_{\mathbb{Y}}^{\ast}v,\eta\rangle_{0,0}$, in both cases we end up with $\int_{0}^{1}\omega(x)\omega'(x)dx=0$, which is a contradiction. Therefore, by \eqref{asyatzero} for the coefficients $b_{j}$, we find that $\underline{c}_{2,j}=0$, $j\in\mathbb{N}_{0}$. We conclude that $\mathrm{Ker}(\lambda-\underline{\Delta}_{\mathbb{Y}}^{\ast})=\{0\}$.

By \cite[Remark 5.26 (i)]{GKM1} the operator $\lambda-\underline{\Delta}_{\mathbb{Y},\min}^{\ast}$ is Fredholm, where $\underline{\Delta}_{\mathbb{Y},\min}^{\ast}$ is the minimal extension of $\Delta_{\mathbb{Y}}$ in $\mathcal{K}_{2}^{0,-\gamma}(\mathbb{Y})$. Since $D(\underline{\Delta}_{\mathbb{Y}}^{\ast})$ differs from $D(\underline{\Delta}_{\mathbb{Y},\min}^{\ast})$ by a finite dimensional space, the operator $\lambda-\underline{\Delta}_{\mathbb{Y}}^{\ast}$ is also Fredholm and hence its range is closed (see e.g. \cite[Proposition 1.3.16]{Lesch}). Since it is also injective, we get that $\lambda-\underline{\Delta}_{\mathbb{Y}}^{\ast}$ is bounded below, which implies that $\bar{\lambda}-\underline{\Delta}_{\mathbb{Y}}$ is surjective. We conclude that $\lambda-\underline{\Delta}_{\mathbb{Y}}$ is bijective, which completes the proof. 
\end{proof}

Based on the above lemma, the holomorphic semigroup generation for $\underline{\Delta}_{\mathbb{Y}}$ is obtained by the following well known argument.

\begin{lemma}\label{sectmc}
Let $p\in(1,\infty)$ and $\gamma$ satisfying \eqref{gammachoice}. Then, for any $\theta\in[0,\pi)$ there exists a $C>1$ such that 
$$
|\lambda|\|(\lambda-\underline{\Delta}_{\mathbb{Y}})^{-1}\|_{\mathcal{L}(\mathcal{K}^{0,\gamma}_{p}(\mathbb{Y}))} \leq C \quad \text{for all} \quad \lambda\in\mathbb{C}\backslash\{0\} \quad \text{with} \quad |\arg(\lambda)|\leq\theta. 
$$
\end{lemma}
\begin{proof} For any $\rho>0$ denote by $\kappa_{\rho}$ the normalized dilation group action on functions on $\mathbb{Y}$ defined by 
$$
(\kappa_{\rho}u)(x,y)=\rho^{\frac{1}{2}-\gamma}u(\rho x,y), \quad x>0.
$$
Then, for any $\lambda\in\mathbb{C}$ and $\rho>0$, in $C_{0}^{\infty}(\mathbb{Y})$ we have
\begin{equation}\label{krhoinv}
\lambda-\Delta_{\mathbb{Y}}=\rho^{2}\kappa_{\rho}(\rho^{-2}\lambda-\Delta_{\mathbb{Y}})\kappa_{\rho}^{-1}.
\end{equation}
Since $D(\underline{\Delta}_{\mathbb{Y}})$ is invariant under $\kappa_{\rho}$, \eqref{krhoinv} still holds in $D(\underline{\Delta}_{\mathbb{Y}})$ if we replace $\Delta_{\mathbb{Y}}$ with $\underline{\Delta}_{\mathbb{Y}}$. Hence, if we restrict $\lambda\in\mathbb{C}\backslash\{0\}$ with $|\arg(\lambda)|\leq\theta$ and choose $\rho=\sqrt{|\lambda|}$, we obtain 
$$
|\lambda|(\lambda-\underline{\Delta}_{\mathbb{Y}})^{-1}=\kappa_{\rho}(|\lambda|^{-1}\lambda-\underline{\Delta}_{\mathbb{Y}})^{-1}\kappa_{\rho}^{-1}
$$
in $\mathcal{K}_{p}^{0,\gamma}(\mathbb{Y})$, where we have used Lemma \ref{nospecminus}. Then, the result follows by the above relation, since $\kappa_{\rho}$ is an isometry on $\mathcal{K}_{p}^{0,\gamma}(\mathbb{Y})$. 
\end{proof}

Finally, the holomorphic semigroup generation for $\underline{\Delta}_{\mathbb{B}}$ is obtained: for $s=0$ by a parametrix construction based on $\underline{\Delta}_{\mathbb{Y}}$, for $s\in \mathbb{N}_{0}$ by induction and for $s\geq0$ by interpolation.

\begin{theorem}\label{sectDB}
Let $p\in(1,\infty)$, $s\geq0$ and $\gamma$ satisfying \eqref{gammachoice}. There exists a $c>0$ such that the spectrum of $\underline{\Delta}_{\mathbb{B}}-c$ is contained in $(-\infty,0]$. Moreover, for any $\theta\in [0,\pi)$, there exists a $C>1$ such that 
\begin{equation}\label{DBsect}
|\lambda|\|(\lambda+c-\underline{\Delta}_{\mathbb{B}})^{-1}\|_{\mathcal{L}(\mathcal{H}^{s,\gamma}_{p}(\mathbb{B}))} \leq C \quad \text{for all} \quad \lambda\in\mathbb{C}\backslash\{0\} \quad \text{with} \quad |\arg(\lambda)|\leq\theta. 
\end{equation}
\end{theorem}
\begin{proof}
Denote $X_{0}^{s}=\mathcal{H}_{p}^{s,\gamma}(\mathbb{B})$, $X_{1}^{s}=\mathcal{H}_{p}^{s+2,\gamma+2}(\mathbb{B})\oplus\mathbb{C}_{\omega}$ and 
$$
\underline{\Delta}_{\mathbb{B}}: X_{1}^{s}\rightarrow X_{0}^{s} \quad \text{by} \quad \underline{\Delta}_{\mathbb{B},s}.
$$
We split the proof into several parts according to different values of $s$.

{\em Case s=0}. Let $\mathbb{F}$ be a closed connected smooth Riemannian manifold such that 
$$
(\mathcal{B}\backslash([0,1/2)\times\partial\mathcal{B}),h|_{\mathcal{B}\backslash([0,1/2)\times\partial\mathcal{B})})
$$ 
is isometrically embedded into $\mathbb{F}$. Denote by $L^{p}(\mathbb{F})$ the space of $p$-integrable functions on $\mathbb{F}$ with respect to the Riemannian measure. Let $\Delta_{\mathbb{F}}$ be the Laplacian on $\mathbb{F}$ and denote by $\underline{\Delta}_{\mathbb{F}}$ the closed extension 
$$
H_{p}^{2}(\mathbb{F})\ni u\mapsto \Delta_{\mathbb{F}}u \in L^{p}(\mathbb{F}). 
$$
The spectrum of $\underline{\Delta}_{\mathbb{F}}$ is contained in $(-\infty,0]$ and, by ellipticity with parameter, see e.g. \cite[Corollary 9.2]{Shu}, there exist a $C_{0}>1$ such that 
\begin{equation}\label{sectDF}
|\lambda|\|(\lambda-\underline{\Delta}_{\mathbb{F}})^{-1}\|_{\mathcal{L}(L^{p}(\mathbb{F}))} \leq C_{0} \quad \text{for all} \quad \lambda\in\mathbb{C}\backslash\{0\} \quad \text{with} \quad |\arg(\lambda)|\leq\theta. 
\end{equation}
Let $\omega_{1}$, $\omega_{2}$ be two cut-off functions with values on $[0,1]$ such that $\omega_{1}=\omega_{2}=1$ on $[0,1/2)\times\partial\mathcal{B}$, $\omega_{1}=\omega_{2}=0$ on $\mathcal{B}\backslash([0,1)\times\partial\mathcal{B})$ and $\omega_{2}=1$ on $\mathrm{supp}(\omega_{1})$. Moreover, let $\omega_{3}=1-\omega_{1}$ and let $\omega_{4}$ be a cut-off function such that $\omega_{4}=0$ on $[0,1/2)\times\partial\mathcal{B}$ and $\omega_{4}=1$ on $\mathrm{supp}(\omega_{3})$. Consider the parametrix 
$$
Q(\lambda)=\omega_{1}(\lambda+c_{1}+c_{0}-\underline{\Delta}_{\mathbb{Y}})^{-1}\omega_{2}+\omega_{3}(\lambda+c_{1}+c_{0}-\underline{\Delta}_{\mathbb{F}})^{-1}\omega_{4}, \quad \lambda\in\mathbb{C}\backslash(-\infty,0),
$$
where $c_{0}, c_{1}>0$ and $\underline{\Delta}_{\mathbb{Y}}$ is the closed extension \eqref{deltahat}. Clearly, $Q(\lambda)$ is a well defined map from $X_{0}^{0}$ to $X_{1}^{0}$. Moreover
\begin{equation}\label{paramtrx}
(\lambda+c_{1}+c_{0}-\underline{\Delta}_{\mathbb{B},0})Q(\lambda)=I-P(\lambda), \quad \lambda\in\mathbb{C}\backslash(-\infty,0),
\end{equation}
where 
\begin{eqnarray}\nonumber
P(\lambda)&=&[\underline{\Delta}_{\mathbb{Y}},\omega_{1}](\lambda+c_{1}+c_{0}-\underline{\Delta}_{\mathbb{Y}})^{-1}\omega_{2}+[\underline{\Delta}_{\mathbb{F}},\omega_{3}](\lambda+c_{1}+c_{0}-\underline{\Delta}_{\mathbb{F}})^{-1}\omega_{4}\\\nonumber
&=&\omega_{2}[\underline{\Delta}_{\mathbb{Y}},\omega_{1}](c_{0}-\underline{\Delta}_{\mathbb{Y}})^{-\eta}(c_{0}-\underline{\Delta}_{\mathbb{Y}})^{\eta}(\lambda+c_{1}+c_{0}-\underline{\Delta}_{\mathbb{Y}})^{-1}\omega_{2}\\\label{gtaet}
&&+\omega_{4}[\underline{\Delta}_{\mathbb{F}},\omega_{3}](c_{0}-\underline{\Delta}_{\mathbb{F}})^{-\eta}(c_{0}-\underline{\Delta}_{\mathbb{F}})^{\eta}(\lambda+c_{1}+c_{0}-\underline{\Delta}_{\mathbb{F}})^{-1}\omega_{4}
\end{eqnarray}
and $\eta\in(1/2,1)$. Here the fractional powers $(c_{0}-\underline{\Delta}_{\mathbb{Y}})^{\xi}$, $(c_{0}-\underline{\Delta}_{\mathbb{F}})^{\xi}$, $\xi\in\mathbb{R}$, are defined by the Dunford functional calculus for sectorial operators, see e.g. \cite[Theorem III.4.6.5]{Am}.

By \cite[(I.2.5.2) and (I.2.9.6)]{Am} and \cite[Lemma 4.5]{RoShao} we have
$$
D((c_{0}-\underline{\Delta}_{\mathbb{Y}})^{\eta})\hookrightarrow [\mathcal{K}^{0,\gamma}_{p}(\mathbb{Y}),\mathcal{K}^{2,\gamma+2}_{p}(\mathbb{Y})\oplus\mathbb{C}_{\omega}]_{\eta-\varepsilon}\hookrightarrow \mathcal{K}^{2\eta-3\varepsilon,\gamma+2\eta-3\varepsilon}_{p}(\mathbb{Y})+\mathbb{C}_{\omega},
$$
for all $\varepsilon>0$ sufficiently small, where $[\cdot,\cdot]_{\eta-\varepsilon}$ stands for complex interpolation. Moreover,
$$
D((c_{0}-\underline{\Delta}_{\mathbb{F}})^{\eta})\hookrightarrow [L^{p}(\mathbb{F}),H_{p}^{2}(\mathbb{F})]_{\eta-\varepsilon}\hookrightarrow H_{p}^{2\eta-3\varepsilon}(\mathbb{F}).
$$
Note that in \eqref{gtaet} the commutator $[\underline{\Delta}_{\mathbb{Y}},\omega_{1}]$ is a first order differential operator on $\mathbb{Y}$ whose coefficients are smooth compactly supported functions that vanish near the boundary $\{0\}\times\partial\mathcal{B}$ and the commutator $[\underline{\Delta}_{\mathbb{F}},\omega_{3}]$ is first order differential operator with smooth coefficients on $\mathbb{F}$. Therefore, the maps 
$$
[\underline{\Delta}_{\mathbb{Y}},\omega_{1}](c_{0}-\underline{\Delta}_{\mathbb{Y}})^{-\eta}: \mathcal{K}^{0,\gamma}_{p}(\mathbb{Y}) \rightarrow \mathcal{K}^{0,\gamma}_{p}(\mathbb{Y})
$$
and
$$
[\underline{\Delta}_{\mathbb{F}},\omega_{3}](c_{0}-\underline{\Delta}_{\mathbb{F}})^{-\eta}:L^{p}(\mathbb{F}) \rightarrow L^{p}(\mathbb{F})
$$
are bounded. Additionally, by \cite[Lemma 7.1]{RoShao} or \cite[Lemma 2.3.3]{Tanabe}, the norms 
$$
\|(c_{0}-\underline{\Delta}_{\mathbb{Y}})^{\eta}(\lambda+c_{1}+c_{0}-\underline{\Delta}_{\mathbb{Y}})^{-1}\|_{\mathcal{L}( \mathcal{K}^{0,\gamma}_{p}(\mathbb{Y}))}, \quad \|(c_{0}-\underline{\Delta}_{\mathbb{F}})^{\eta}(\lambda+c_{1}+c_{0}-\underline{\Delta}_{\mathbb{F}})^{-1}\|_{\mathcal{L}(L^{p}(\mathbb{F}))},
$$
become arbitrary small, uniformly in $\lambda$, by taking $c_{1}$ sufficiently large. Therefore, by \eqref{paramtrx}-\eqref{gtaet} we obtain a right inverse of $\lambda+c_{1}+c_{0}-\underline{\Delta}_{\mathbb{B},0}$ which belongs to $\mathcal{L}(X_{0}^{0},X_{1}^{0})$. 

Similarly, on $X_{1}^{0}$ we have
\begin{equation}\label{leftinvpar}
Q(\lambda)(\lambda+c_{1}+c_{0}-\underline{\Delta}_{\mathbb{B},0})=I-R(\lambda), \quad \lambda\in\mathbb{C}\backslash(-\infty,0),
\end{equation}
where
\begin{eqnarray}\nonumber
R(\lambda)&=&\omega_{1}(\lambda+c_{1}+c_{0}-\underline{\Delta}_{\mathbb{Y}})^{-1}[\omega_{2},\underline{\Delta}_{\mathbb{Y}}]+\omega_{3}(\lambda+c_{1}+c_{0}-\underline{\Delta}_{\mathbb{F}})^{-1}[\omega_{4},\underline{\Delta}_{\mathbb{F}}]\\\nonumber
&=&\omega_{1}(\lambda+c_{1}+c_{0}-\underline{\Delta}_{\mathbb{Y}})^{-1} (c_{0}-\underline{\Delta}_{\mathbb{Y}})^{-\rho}(c_{0}-\underline{\Delta}_{\mathbb{Y}})^{\rho}[\omega_{2},\underline{\Delta}_{\mathbb{Y}}]\\\label{Rterm}
&&+\omega_{3}(\lambda+c_{1}+c_{0}-\underline{\Delta}_{\mathbb{F}})^{-1}(c_{0}-\underline{\Delta}_{\mathbb{F}})^{-\rho}(c_{0}-\underline{\Delta}_{\mathbb{F}})^{\rho}[\omega_{4},\underline{\Delta}_{\mathbb{F}}]
\end{eqnarray}
and $\rho\in(0,1/2)$. 
By \cite[(I.2.5.2) and (I.2.9.6)]{Am} and \cite[Lemma 4.5]{RoShao} we have
$$
 \mathcal{K}^{2\rho+3\varepsilon,\gamma+2\rho+3\varepsilon}_{p}(\mathbb{Y})+\mathbb{C}_{\omega}\hookrightarrow [\mathcal{K}^{0,\gamma}_{p}(\mathbb{Y}),\mathcal{K}^{2,\gamma+2}_{p}(\mathbb{Y})\oplus\mathbb{C}_{\omega}]_{\rho+\varepsilon}\hookrightarrow D((c_{0}-\underline{\Delta}_{\mathbb{Y}})^{\rho}),
$$
for all $\varepsilon>0$ small enough. Moreover,
$$
H_{p}^{2\rho+3\varepsilon}(\mathbb{F}) \hookrightarrow [L^{p}(\mathbb{F}),H_{p}^{2}(\mathbb{F})]_{\rho+\varepsilon} \hookrightarrow D((c_{0}-\underline{\Delta}_{\mathbb{F}})^{\rho}).
$$
Hence, the maps 
$$
(c_{0}-\underline{\Delta}_{\mathbb{Y}})^{\rho}[\omega_{2},\underline{\Delta}_{\mathbb{Y}}]: \mathcal{K}^{2,\gamma+2}_{p}(\mathbb{Y})\oplus\mathbb{C}_{\omega} \rightarrow \mathcal{K}^{0,\gamma}_{p}(\mathbb{Y})
$$
and
$$
(c_{0}-\underline{\Delta}_{\mathbb{F}})^{\rho}[\omega_{4},\underline{\Delta}_{\mathbb{F}}]: H_{p}^{2}(\mathbb{F}) \rightarrow L^{p}(\mathbb{F})
$$
are bounded. Moreover, by \cite[Lemma 7.1]{RoShao} or \cite[Lemma 2.3.3]{Tanabe}, the norms
$$
\|(\lambda+c_{1}+c_{0}-\underline{\Delta}_{\mathbb{Y}})^{-1} (c_{0}-\underline{\Delta}_{\mathbb{Y}})^{-\rho}\|_{\mathcal{L}(\mathcal{K}^{0,\gamma}_{p}(\mathbb{Y}),\mathcal{K}^{2,\gamma+2}_{p}(\mathbb{Y})\oplus\mathbb{C}_{\omega})}
$$
and
$$
\|(\lambda+c_{1}+c_{0}-\underline{\Delta}_{\mathbb{F}})^{-1}(c_{0}-\underline{\Delta}_{\mathbb{F}})^{-\rho}\|_{\mathcal{L}(L^{p}(\mathbb{F}),H_{p}^{2}(\mathbb{F}))}
$$
become arbitrary small, uniformly in $\lambda$, by taking $c_{1}$ sufficiently large. As a consequence, by \eqref{leftinvpar}-\eqref{Rterm} we obtain a left inverse of $\lambda+c_{1}+c_{0}-\underline{\Delta}_{\mathbb{B},0}$ that belongs to $\mathcal{L}(X_{0}^{0},X_{1}^{0})$. 

We conclude that
$$
(\lambda+c-\underline{\Delta}_{\mathbb{B},0})^{-1}=Q(\lambda)(I-P(\lambda))^{-1}\in \mathcal{L}(X_{0}^{0},X_{1}^{0}), \quad \lambda\in \mathbb{C}\backslash(-\infty,0),
$$
where $c=c_{1}+c_{0}$. The result for $s=0$ follows by Lemma \ref{sectmc}, \eqref{sectDF} and the above formula. 

{\em Case $s\in\mathbb{N}_{0}$}. First we show that, for any $\nu>0$ we have that
\begin{equation}\label{resolvrstr}
\mathbb{C}\backslash(-\infty,0)\in \rho(\underline{\Delta}_{\mathbb{B},\nu}-c) \quad \text{and} \quad (\lambda+c-\underline{\Delta}_{\mathbb{B},\nu})^{-1}= (\lambda+c-\underline{\Delta}_{\mathbb{B},0})^{-1}|_{X_{0}^{s}}, \quad \lambda\in \mathbb{C}\backslash(-\infty,0).
\end{equation}
It suffices to verify that the identities
$$
(\lambda+c-\underline{\Delta}_{\mathbb{B},0})^{-1}(\lambda+c-\underline{\Delta}_{\mathbb{B},\nu})=I \quad \text{and} \quad (\lambda+c-\underline{\Delta}_{\mathbb{B},\nu})(\lambda+c-\underline{\Delta}_{\mathbb{B},0})^{-1}=I
$$
hold on $X_{1}^{\nu}$ and $X_{0}^{\nu}$ respectively. The first one is trivial. For the second one, restrict $\nu\in(0,2]$ and let $u\in X_{1}^{0}$ such that $(\lambda+c-\Delta_{\mathbb{B}})u\in X_{0}^{\nu}$. Then
$$
u\in D(\underline{\Delta}_{\mathbb{B},\nu,\max}) \cap X_{1}^{0},
$$
where $D(\underline{\Delta}_{\mathbb{B},\nu,\max})$ is the maximal domain of $\Delta_{\mathbb{B}}$ in $X_{0}^{\nu}$. By the structure of $D(\underline{\Delta}_{\mathbb{B},\nu,\max})$ we infer that $u\in X_{1}^{\nu}$. The result then follows by iterating the above procedure. 

Assume that \eqref{DBsect} holds for some $s\in\mathbb{N}_{0}$. We proceed by induction in $s$. Let $U_{j}$, $j\in\{1,\dots,\ell\}$, $\ell\in\mathbb{N}$, be a covering of $\mathcal{B}$ by coordinate charts and let $\{\phi_{j}\}_{j\in\{1,\dots,\ell\}}$ be a subordinate partition of unity, where we choose $U_{1}=[0,r)\times\partial\mathcal{B}$, for certain $r\in(0,1)$. Let $(z_{1},\dots,z_{n+1})$ be local coordinates in $U_{j}$ such that, when $j=1$ we choose $z_{1}=x$, $(z_{2},\dots,z_{n+1})=y\in \partial\mathcal{B}$ and by $\partial_{z_{1}}$ we denote $x\partial_{x}$. For any $v\in X_{0}^{s+1}$ we have
\begin{eqnarray*}
\lefteqn{\|(\lambda+c-\underline{\Delta}_{\mathbb{B},s+1})^{-1}v\|_{X_{0}^{s+1}}}\\
&=&\|(\lambda+c-\underline{\Delta}_{\mathbb{B},s})^{-1}v\|_{X_{0}^{s+1}}\\
&\leq&\sum_{j=1}^{\ell}\sum_{k=1}^{n+1}\sum_{|\alpha|\leq1}\|\partial_{z_{k}}^{\alpha}\phi_{j}(\lambda+c-\underline{\Delta}_{\mathbb{B},s})^{-1}v\|_{X_{0}^{s}}\\
&\leq&\sum_{j=1}^{\ell}\sum_{k=1}^{n+1}\sum_{|\alpha|\leq1}\Big( \|(\lambda+c-\underline{\Delta}_{\mathbb{B},s})^{-1}\partial_{z_{k}}^{\alpha}\phi_{j}v\|_{X_{0}^{s}}+\|[\partial_{z_{k}}^{\alpha}\phi_{j},(\lambda+c-\underline{\Delta}_{\mathbb{B},s})^{-1}]v\|_{X_{0}^{s}}\Big)\\
&=&\sum_{j=1}^{\ell}\sum_{k=1}^{n+1}\sum_{|\alpha|\leq1}\Big( \|(\lambda+c-\underline{\Delta}_{\mathbb{B},s})^{-1}\partial_{z_{k}}^{\alpha}\phi_{j}v\|_{X_{0}^{s}}\\
&&+\|(\lambda+c-\underline{\Delta}_{\mathbb{B},s})^{-1}[\underline{\Delta}_{\mathbb{B},s},\partial_{z_{k}}^{\alpha}\phi_{j}](c-\underline{\Delta}_{\mathbb{B},s})^{-1}(c-\underline{\Delta}_{\mathbb{B},s})(\lambda+c-\underline{\Delta}_{\mathbb{B},s})^{-1}v\|_{X_{0}^{s}}\Big).
\end{eqnarray*}
Note that $[\underline{\Delta}_{\mathbb{B},s},\partial_{z_{k}}^{\alpha}\phi_{j}]$ are second order cone differential operators whose zero order coefficients have support away of $\{0\}\times\partial\mathcal{B}$. Therefore $[\underline{\Delta}_{\mathbb{B},s},\partial_{z_{k}}^{\alpha}\phi_{j}](c-\underline{\Delta}_{s})^{-1}$ belong to $\mathcal{L}(X_{0}^{s})$. Consequently, we have
\begin{eqnarray*}
\lefteqn{\|(\lambda+c-\underline{\Delta}_{\mathbb{B},s+1})^{-1}v\|_{X_{0}^{s+1}}}\\
&\leq&\|(\lambda+c-\underline{\Delta}_{\mathbb{B},s})^{-1}\|_{\mathcal{L}(X_{0}^{s})}\Big[\|v\|_{X_{0}^{s+1}}+2\ell(n+1)\Big(\max_{j,k,\alpha}\Big\{\|[\underline{\Delta}_{\mathbb{B},s},\partial_{z_{k}}^{\alpha}\phi_{j}](c-\underline{\Delta}_{\mathbb{B},s})^{-1}\|_{\mathcal{L}(X_{0}^{s})}\Big\}\Big)\\
&&\times\|(c-\underline{\Delta}_{\mathbb{B},s})(\lambda+c-\underline{\Delta}_{\mathbb{B},s})^{-1}\|_{\mathcal{L}(X_{0}^{s})}\|v\|_{X_{0}^{s}}\Big],
\end{eqnarray*}
which shows the result for $s+1$.

{\em Case $s\geq0$}. By \cite[Lemma 2.12]{Roidos2} and \cite[Lemma 3.7]{RS1}, for any $\sigma\in[0,1]$ we have $[X_{0}^{s},X_{0}^{s+1}]_{\sigma}=X_{0}^{s+\sigma}$ and $[X_{1}^{s},X_{1}^{s+1}]_{\sigma}=X_{1}^{s+\sigma}$. Then the result follows by \eqref{resolvrstr} and interpolation, see e.g. \cite[Theorem 2.6]{Lunardi2}.
\end{proof}

\begin{remark}
We can extend Theorem \ref{sectDB} to the case of $s\in\mathbb{R}$. This can be achieved by applying the above steps to the adjoint $\underline{\Delta}_{\mathbb{B},-s}^{\ast}:D(\underline{\Delta}_{\mathbb{B},-s}^{\ast})\rightarrow \mathcal{H}_{q}^{-s,-\gamma}(\mathbb{B})$ of $\underline{\Delta}_{\mathbb{B},s}$, $p^{-1}+q^{-1}=1$, $s>0$, and then using duality, see e.g. \cite[Proposition 1.3 (v)]{DHP}.
\end{remark}

\end{document}